\documentclass[11pt]{amsart}
\usepackage{amssymb}
\usepackage{palatino}
\usepackage{tikz-cd}
\input amssym.def

\usepackage{xcolor}
\usepackage{amsmath, amsfonts}
\usepackage{amssymb}
\usepackage{amscd}
\usepackage[mathscr]{eucal}
\usepackage{palatino}
\usepackage{theoremref}
\usepackage{tikz}
\usepackage{float}
\usepackage{eqnarray}

\newfont{\cyrr}{wncyr10}

\newcommand{\thmref}[1]{Theorem~\ref{#1}}

\newcommand{\propref}[1]{Proposition~\ref{#1}}

\newtheorem{thm}{Theorem}

\newtheorem{prop}[thm]{Proposition}
\newtheorem{rmk}{Remark}[section]

\newtheorem*{conj}{Conjecture}

\newcommand{\Z}{{\mathbb Z}}

\def\({\left(}
\def\){\right)}
\def\[{\left[}
\def\]{\right]}

\def\N{\mathbb{N}}
\def\R{\mathbb{R}}
\def\C{\mathbb{C}}

\def\cA{\mathcal{A}}

\def\cE{\mathcal{E}}
\def\cF{\mathcal{F}}
\def\cG{\mathcal{G}}

\def\cP{\mathcal{P}}

\def\cS{\mathcal{S}}
\def\cT{\mathcal{T}}
\def\cG{\mathcal{G}}

\def\cS{\mathcal{S}}

\def\e{\epsilon}

\title{Brun's inequality for a geometric lattice}
\author{M. Ram Murty and Sunil Naik}

\address{ M. Ram Murty and Sunil Naik \newline
	Department of Mathematics,
	Queen's University, Jeffery Hall, 
	99 University Avenue, 
	Kingston, ON K7L 3N6, 
	Canada}
\email{murty@queensu.ca}
\email{naik.s@queensu.ca}

\begin{document}
	
\hfuzz 5pt	


	
\maketitle

\begin{abstract}
%We generalize Brun's inequality of classical sieve theory in the context of a geometric lattice.  We then apply this to study
%Dowling lattices.
In a seminal paper of 1915, V. Brun introduced Brun's sieve, which is based on Brun's inequality for the M\"{o}bius function and is a very powerful tool in modern number theory. The importance of the M\"{o}bius function in enumeration problems led G.-C. Rota to introduce the concept of the M\"{o}bius function to partially ordered sets. In this article, we prove Brun's inequality for geometric lattices and develop a sieve in this context. One of the main ingredients is a recent work of K. Adiprasito, J. Huh, and E. Katz on the log-concavity of absolute values of the Whitney numbers associated with matroids. We also study shifted convolutions of the Whitney numbers associated with Dowling lattices. Further, we derive an asymptotic formula for generalized Dowling numbers.
\end{abstract}

\section{Introduction}

One of the most powerful tools in analytic number theory is
the sieve method.  Classically, the sieve methods of Eratosthenes,
Brun and Selberg are combinatorial in nature, whereas the large
sieve inequality due to Linnik is Fourier analytic.  The former class
of sieve methods could, in principle, be formulated in a general 
combinatorial setting.  The first step in this direction was
due to Wilson \cite{wilson1}, \cite{wilson2}
who formulated a  ``Selberg sieve inequality''
for partially ordered sets satisfying some standard properties.
Later, Liu and Murty
\cite{liu-murty}
isolated an idea of Tur\'an and formulated a very general
combinatorial sieve that they called ``the Tur\'an sieve'' and 
applied it to an assortment of questions in graph theory and combinatorics.  Though interesting, these results have had 
limited impact.  It would seem that the classical Brun's sieve
would afford a similar generalization.  However, this was not
the case since the simple seed idea which generates the classical Brun's sieve
required a unimodal property of certain Whitney numbers before it could
be generalized to a geometric lattice.  The required unimodal property
is now a theorem due to
June Huh and his school.  Using this theory, we generalize Brun's
inequalities and then develop a general combinatorial sieve.
We  apply this to the study of Dowling lattices.

In 1919, Brun developed his sieve method by making the following
elementary observation concerning the classical M\"obius function
$\mu(n)$.  Let $\omega(n)=r$ denote the number of distinct prime factors of $n$.  Then,
$$ \sum_{d|n \atop \omega(d)\leq k} \mu(d) = (-1)^k {r-1\choose k}$$
so the sum on the left is positive when $k$ is even and negative
when $k$ is odd.  The sieve method is then developed from here.

To generalize this to a combinatorial setting, one would require
some form of ``rank'' function that would serve as an analog of
$\omega(n)$.  The other difficulty is to probe an oscillatory theorem
similar to the one stated for the classical case.  This is the content
of our first theorem described in the next section.

\section{Statements of Results}

The concept of a M\"obius function for a partially ordered set
originates in the 1964 work of Rota \cite{Ro64}.  We refer the reader
to this paper for the basic background though we give a brief
review in later sections.

\begin{thm}[Brun's inequality for a geometric lattice]\label{ThmB}
Let $L$ be a geometric lattice. Then  for any positive integer $k$, we have
$$
\sum_{\substack{x \in L \\ r(x) \leq 2k}} \mu(\hat{0}, x) 
~\geq~ 0
\phantom{mm} \text{and} \phantom{mm}
\sum_{\substack{x \in L \\ r(x) \leq 2k-1}} \mu(\hat{0}, x) 
~\leq~ 0.
$$
Here $\mu$ denotes the M\"obius function on $L$ and $r$ denotes a rank function on $L$.
\end{thm}

The following theorem on unimodal sequences is an important tool for incorporating the work of Adiprasito-Huh-Katz in our context.
\begin{thm}\label{ThmU}
Let $\{a_i\}_{i=0}^{n}$ be a unimodal sequence of 
non-negative real numbers such that
$$
\sum_{i=0}^{n} (-1)^i a_i ~=~ 0.
$$
Then for any positive integer $k$, we have
$$
\sum_{i=0}^{2k} (-1)^i a_i 
~\geq~ 0
\phantom{mm} \text{and} \phantom{mm}
\sum_{i=0}^{2k-1} (-1)^i a_i 
~\leq~ 0.
$$
\end{thm}

\medspace

We use the properties of the M\'{o}bius function of partially ordered sets to develop the following combinatorial sieve for geometric lattices.
\begin{thm}\label{thmSift}
Let $L$ be a geometric lattice of rank $n$ and $\cA \subseteq L$. 
Also let $\cT$ be a set of atoms in $L$ whose join is $\tau$. 
Suppose there exists a function
$$
f ~:~ \N \cup\{0\}  ~\longrightarrow~ [0, \infty)
$$
and a positive real number $X$ such that
\begin{equation*}
\#\{ a \in \cA ~:~ a \geq y\} 
~=~ f(\text{cr}(y)) X ~+~ \cE(y)
\end{equation*}
and 
\begin{equation*}
\cE(y) ~\ll~ \text{cr}(y) f(\text{cr}(y)),
\end{equation*}
where $\text{cr}(y)$ denotes the co-rank of $y$ 
given by $\text{cr}(y) = n- r(y)$. Then we have
\begin{equation*}
\cS(\cA, \cT) 
~=~ \# \{a \in \cA ~:~ a \wedge \tau ~=~ \hat{0}\} 
~=~ X \sum_{k=0}^{\text{r}(\tau)} f(n-k) w_{k}\([\hat{0}, \tau]\) ~+~ E(\cA, \cT),
\end{equation*}
	where
$$
E(\cA, \cT) ~\ll~ 
\sum_{k=0}^{\text{r}(\tau)} (n-k)f(n-k) |w_{k}([\hat{0}, \tau])|.
$$
\end{thm}

\begin{rmk}
One can apply \thmref{ThmB} to derive lower and upper bounds for $\cS(\cA, \cT)$.
\end{rmk}

\medspace

Now we apply our results to investigate Dowling lattices. This leads us to the study of shifted convolution of Whitney numbers associated with these lattices. Let $n$ be a positive integer and $G$ be a multiplicative group of order $m$. 
Let $Q_n(G)$ be the Dowling lattice of rank $n$ (see subsection \ref{SSDowling}). 
It is well-known that $Q_n(G)$ is  a geometric lattice.  
We study relations between certain sifted sets  
and shifted convolution of Whitney numbers of these lattices. Let
$$
w_m(n,k) ~=~ \sum_{x \in Q_n(G) \atop r(x) = n-k} \mu(\hat{0}, x)
\phantom{m}\text{and}\phantom{m} W_m(n,k) 
~=~ \sum_{x \in Q_n(G) \atop r(x) = n-k} 1,
$$ 
where  $w_m(n,k)$ and $W_m(n,k)$ denote the Whitney numbers of $Q_n(G)$ 
of the first kind and the second kind respectively. 
Here $\mu$ denotes the M\"obius function of $Q_n(G)$ 
and $r$ denotes a rank function on $Q_n(G)$. 
In this context, we prove the following result.

\begin{prop}\label{propshiftW}
Let $n$ be a positive integer. For non-negative integers $s, t$, let
\begin{equation*}
c_{n,t}(s) ~=~ \sum_{k \geq 0} w_m(n, k) W_m(k+s,t),
\end{equation*}	
denote shifted convolution of Whitney numbers of $Q_n(G)$. Then we have
$$
c_{n, t}(s) ~=~ 0 \phantom{m}\text{if}\phantom{m} t < n.
$$
If $t \geq n$, then we have
$$
\sum_{s \geq 0} c_{n,t}(s) x^s 
~=~ \frac{1}{x} \prod_{j=n}^{t} \frac{x}{1-(1+jm)x}.
$$
In particular, we have
$$
\sum_{k \geq 0} w_m(n, k) W_m(k,t) ~=~ \delta(n, t).
$$
Here $\delta(n, t)=1$ if $n=t$ and $\delta(n, t)=0$ otherwise.
\end{prop}

\medspace

Further, the study of estimation of sifted sets and 
shifted convolution of Whitney numbers of $Q_n(G)$ 
leads us to the study of $r$-Dowling numbers $D_{m,r}(n)$ 
(see  section \ref{SDmr}). 
In this setup, we prove the following asymptotic formula for $r$-Dowling numbers.
\begin{thm}\label{thmexpDn}
Let $m, n$ and $r$ be positive integers. For $n > e^m$, we have
$$
D_{m,r}(n) ~=~ 
\frac{e^{g_0}}{\sqrt{4 \pi g_2}} \cdot \frac{n!}{\delta^n} 
\(1 ~+~ O\( \frac{(\log n)^{12}}{\sqrt{n}}\)\),
$$
where $\delta$ is a positive real number such that
$
\delta ( r + e^{m \delta}) ~=~ n
$
and  the implied constant is absolute.
Here		
$$
g_0 ~=~ r \delta  ~+~ \frac{e^{m\delta} -1}{m}~
\phantom{mm}\text{and}\phantom{mm}
g_2 ~=~ \frac{n ~+~ m \delta^2 e^{m \delta}}{2}.
$$
	\end{thm}

\medspace

\section{Preliminaries}

\subsection{Unimodality and log-concavity}
A sequence $\{a_i\}_{i=0}^{n}$ of real numbers 
is said to be {\bf unimodal} if there exists an index 
$0 \leq j \leq n$ such that
$$
a_0 ~\leq~ a_1 ~\leq~ \cdots ~\leq~ a_{j-1} 
~\leq~ a_j ~\geq~ a_{j+1} ~\geq~ \cdots ~\geq~ a_n.
$$
There are numerous naturally occurring sequences 
which are unimodal. For example, the sequence 
$\{\binom{n}{k}\}_{k=0}^{n}$ of binomial coefficients is unimodal.

A sequence $\{a_i\}_{i=0}^{n}$ is said to be {\bf log-concave} 
if
$$
a_k^2 ~\geq~ a_{k-1} a_{k+1} 
~\phantom{mm}
\text{for all } ~ 0 < k < n.
$$
The sequence of binomial coefficients is log-concave. 
One can easily show that any log-concave sequence is unimodal. 
Other important examples of log-concave sequences 
that occur in combinatorics are 
{\bf Stirling numbers of the first kind}. 
Stirling numbers $s(n,k)$ of the first kind are defined by
$$
x(x-1) \cdots (x-(n-1)) ~=~ \sum_{k=1}^{n} s(n, k) x^k.
$$
Note that $|s(n,k)|$ is the number of permutations 
of the symmetric group $S_n$ which can be written 
as a product of $k$-disjoint cycles. 
One can see the log-concavity of these numbers 
from an important result of Isaac Newton (see \cite[p. 52]{HLP} \cite{Ne}) 
which states that if all the roots of a polynomial 
$$
f(x) ~=~ \sum_{k=0}^{n} a_k x^k
$$
are real, then the sequence $\{a_k\}_{k=0}^{n}$ 
of coefficients of $f$ is log-concave.
This result produces a rich number of log-concave sequences, 
since the characteristic polynomial 
of a real symmetric matrix has all its roots real.

\medspace

\subsection{Prerequisites from theory of partially ordered sets (posets)}
Let $(P, \leq)$ be a {\bf partially ordered set} ({\bf poset}). 
For any two elements $x, z \in P$, the {\bf interval} $[x, z]$ is given by
$$
[x, z] ~=~ \{ y \in P ~:~ x \leq y \leq z\}.
$$
 A poset $P$ is called {\bf locally finite} 
 if every interval of $P$ is a finite set. 
 In 1964, Rota \cite{Ro64} introduced the concept of 
 M\"obius function on locally finite posets and 
 proved a fundamental property of the M\"obius function, 
 namely M\"obius inversion formula and 
 thereby initiating a combinatorial theory 
 of locally finite posets . 
 The {\bf M\"obius function} $\mu : P \times P \rightarrow \Z$ 
 is defined recursively as follows: 
 $\mu(x, z) = 0$ if $x \not\leq z$, 
 $\mu(x, x)=1$ for all $x \in P$ and 
 if $x<z$,
 $$
 \mu(x, z) ~=~ - \sum_{x \leq y < z} \mu(x, y).
 $$
It is clear from the definition that
$$
\sum_{x \leq y \leq z} \mu(x, y) 
~=~ \delta(x, z),
$$
where $\delta$ denotes the {\bf Kronecker delta function} 
given by
\begin{equation*}
\delta(x, z) 
~=~ 
\begin{cases}
	1 & \text{if  ~} x = z,\\
	0 & \text{otherwise}.
\end{cases}
\end{equation*}

Given two  elements $x, y \in P$, 
we say that $z \in P$ is an {\bf upper bound} for $x$ and $y$ 
if $x \leq z$ and $y \leq z$. 
We say that $z$ is a {\bf least upper bound} for $x$ and $y$ 
if $z$ is an upper bound for $x$ and $y$ and 
if $w$ is an upper bound for $x$ and $y$, then $z \leq w$. 
Clearly, if a least upper bound  for $x$ and $y$ exists, 
then it is unique and is denoted by $x \vee y$ 
(called  $x$ {\bf join} $y$). 
In a similar way, one defines the {\bf greatest lower bound} 
for $x$ and $y$ (also called {\bf meet} of $x$ and $y$) and 
it is denoted by $x \wedge y$, if it exists. 
A poset $P$ is said to have a {\bf least element} 
(usually denoted by $\hat{0}$) if $\hat{0} \leq x$ for all $x \in P$. 
In a similar way, a poset $P$ is said to have a {\bf greatest element} 
(usually denoted by $\hat{1}$) 
if $x \leq \hat{1}$ for all $x \in P$. 

A {\bf lattice} $L$ is a poset in which the least upper bound 
and the greatest lower bound for any two elements of $L$ exist. 
Now onwards, we will focus on finite lattices. 
Notice that a finite lattice has least and greatest elements. 
An element $\hat{0} \neq x \in L$ is called an {\bf atom} 
if there exists no element $y \in L$ with $\hat{0} < y < x$, 
in other words $[\hat{0}, x] = \{\hat{0}, x\}$. 
A lattice is said to be {\bf atomistic} if every non-zero element 
is a join of atoms. Given $x < y$ in $L$, 
we say that $y$ {\bf covers} $x$ if $[x, y] = \{x, y\}$ and 
we denote it by $x <: y$. 
A lattice $L$ is said to be {\bf graded} 
if there exists a function
$$
r ~:~ L ~\longrightarrow~ \R_{\geq 0}
$$
such that $r(\hat{0}) = 0$ and 
$r(y) = r(x)+1$ if $x <: y$. 
A graded lattice is said to be {\bf semimodular} if 
$$
r(x \wedge y) + r(x \vee y) ~\leq~ r(x) + r(y)
$$
for any $x, y \in L$. A {\bf geometric lattice} 
is a finite atomistic, semimodular lattice. 
For example, the set $\Pi_n$ of all partitions of $\{1, 2, \cdots, n\}$ 
ordered by refinement and the set of all subspaces of 
an $n$-dimensional vector space over a finite field 
ordered by inclusion are geometric lattices 
(see \cite[Ch. 6]{CM} \cite{RM} \cite[Ch. 3]{St} for more details).

\medspace

\subsection{Prerequisites from theory of matroids}
A {\bf matroid} $M$ is pair $\(E, \cF\)$, 
where $E$ is a finite set and $\cF$ is a collection of subsets of $E$ 
(called {\bf independent sets}) 
satisfying the following axioms:
\begin{enumerate}
	\item[i)] $\emptyset \in \cF$;
	
	\item[ii)] If $A \in \cF$ and $B \subseteq A$, then $B \in \cF$;
	
	\item[iii)] If $A , B \in \cF$ and $|A| < |B|$, 
	then there exists $b \in B \backslash A$ 
	such that $A \cup \{b\} \in \cF$ ({\bf exchange property}).
\end{enumerate}

A maximal independent set of $E$ is called a {\bf basis} of $M$. 
As a consequence of exchange property, 
all bases of $M$ have the same cardinality 
and is called the {\bf rank} of $M$ and is denoted by $r(M)$. 
For any $A \subseteq E$, the rank of $A$ is defined to be 
the cardinality of a maximal independent subset of $A$ and is denoted by $r(A)$. 
A subset of $E$ which is not independent is called 
a {\bf dependent set} and 
a minimal dependent set is called a {\bf circuit}.

One defines the {\bf characteristic polynomial} 
of a matroid $M$ as follows:
$$
\chi_{M}(\lambda) ~=~
 \sum_{A \subseteq E } (-1)^{|A|} \lambda^{r(M)-r(A)}.
$$
As the name suggests, $\chi_{M}(\lambda)$ is a polynomial 
in $\lambda$ of degree equal to the rank of $M$ 
and all of its coefficients are integers. 
Hence one can write it as
$$
\chi_{M}(\lambda) ~=~ \sum_{i=0}^{r(M)} w_i(M) \lambda^{r(M) - i}.
$$
 The above integers $w_i(M)$ are called 
 {\bf Whitney numbers of the first kind} associated with the matroid $M$. 
 Rota \cite{Ro71} conjectured that the sequence of absolute values 
 of these numbers is {\bf unimodal}. 
 A stronger conjecture was proposed by Heron \cite{He} and  Welsh \cite{Wel}.
 \begin{conj}[Heron-Rota-Welsh]
 The sequence $\{|w_i(M)|\}_{i=0}^{r(M)}$ of absolute values of the coefficients 
 of $\chi_M(\lambda)$ is log-concave.
 \end{conj}
 
This conjecture was resolved in 2018 by Adiprasito-Huh-Katz \cite{AHK}, 
one of the works for which J. Huh has been awarded Fields medal in 2022.

\vspace{2mm}

We have seen the definition of a matroid using independent axioms 
which abstract the notion of linear independence. 
Now we will see an another way to define a matroid
in terms of {\bf closure operators} 
which abstract the notion of linear span. 
Let $\cP(E)$ be the power set of $E$. 
A matroid $M$ is a finite set $E$ together with a function
\begin{equation*}
	\begin{split}
\text{cl} ~:~ \cP(E) ~&\longrightarrow~ \cP(E) \\
                X ~&\longmapsto~ \overline{X}
            \end{split}
\end{equation*}
satisfying the following axioms:
\begin{enumerate}
	\item[i)] $X \subseteq \overline{X}$;
	
	\item[ii)] If $Y \subseteq X$, 
	then $\overline{Y} \subseteq \overline{X}$;
	
	\item[iii)] $\overline{\overline{X}} =  \overline{X}$;
	
	\item[iv)] If $y \in \overline{X \cup \{x\}}$ 
	and $y \not\in \overline{X}$, 
	then $x \in \overline{X \cup \{y\}}$.
\end{enumerate}
The above two definitions of a matroid are equivalent. 
To see this, given a matroid with independent axioms, 
one defines the closure of a set  $X \subseteq E$ by
$$
\overline{X} ~:=~ \{x \in E : r(X \cup \{x\}) = r(X)\}.
$$ 
It is easy to check that the closure operator satisfies the span axioms. 
Conversely, given a matroid with span axioms, 
one declares a subset $A \subseteq E$ to be independent 
if $x \in A$ implies $x \not\in \overline{A\backslash\{x\}}$.

A subset $X \subseteq E$ is called {\bf flat} if $\overline{X} = X$. 
The elements of $\overline{\emptyset}$  are called {\bf loops}. 
Two elements $x, y \in E$ are said to be {\bf parallel} 
if $\{x, y\}$ is a circuit. 
A matroid is said to be {\bf simple} if it has no loops and parallel points. 
Every matroid has a canonical {\bf simplification} $\hat{M}$ 
obtained by removing loops and identifying parallel points. 
More precisely,
if $E'$ denotes the set of all elements of $E$ which are not loops, 
then one defines an equivalence relation on $E'$ by:
$$
 x ~\sim~ y ~~\iff~~ \text{$x$ and $y$ are parallel}.
$$
Let $\hat{E}$ denotes the set of all equivalence classes of $E'$, 
then $\hat{E} = \{\overline{x} : x \in E'\}$.
We say a subset $\{\overline{x_1}, \overline{x_2}, 
\cdots, \overline{x_k}\} \subseteq \hat{E}$ 
is independent if $\{x_1, x_2, \cdots, x_k\}$ is independent. 
Let $\hat{\cF}$ be the set of all independent subsets of $\hat{E}$ 
along with empty set, then one can check that 
$\hat{M} = (\hat{E}, \hat{\cF})$ is a simple matroid. 
Please refer to \cite{Ba, RM, Ox} for more details.

\medspace

\subsection{Prerequisites from Dowling lattices}\label{SSDowling}
Let $X=\{x_1, x_2, \cdots , x_n\}$ be a finite set of $n$ elements. 
Let $P_n$ be the collection of partitions of $X$. 
Then the set $P_n$ is partially ordered by {\bf refinement} : 
we say $\alpha \leq \beta$, if every block of $\beta$ 
is a union of blocks of $\alpha$. For example,
$$
\{x_1,x_2\}\{x_3\}\{x_4\}\{x_5, x_6, \cdots, x_n\} 
~\leq~ \{x_1,x_2,x_3\}\{x_4\}\{x_5, x_6, \cdots, x_n\}.
$$ 
By a {\bf partial partition} of $X$, we mean a collection 
$\cF= \{A_1, A_2,\cdots, A_r\}$ of mutually disjoint non-empty subsets of $X$. 
For example, if we take $A_1= \{x_1, x_2\}$ and $A_2 = \{x_n\}$, 
then $\{A_1, A_2\}$ is a partial partition of $X$. 
The subsets $A_i$'s are called {\bf blocks} of $\cF$. 
Let $Q_n$ be the set of all partial partitions of $X$. 
As before, the set $Q_n$ is partially ordered by refinement. 
Let  $x_0$ be an element  which does not belong to $X$.
We have an isomorphism
\begin{align*}
 Q_n &~\longrightarrow~ P_{n+1} \\
      \{A_1, A_2, \cdots, A_r\} 
      &~\longmapsto~ \{A_0 \cup \{x_0\}, A_1, A_2, \cdots, A_r\},
\end{align*}
where $A_0 = X \backslash \(\cup_{i=1}^{r} A_i\)$. 
The block $A_0 \cup \{x_0\}$ is called the {\bf zero block} 
the partition $\{A_0 \cup \{x_0\}, A_1, A_2, \cdots, A_r\}$ 
of $X \cup \{x_0\}$. If $\cF$ is a partial partition of $X$, 
then the rank of $\cF$ is given by
$$
r(\cF) ~=~ n - |\cF|,
$$
where $|\cF|$ denotes the number of blocks of $\cF$.

Let $G$ be a finite multiplicative group. 
By a {\bf partial $G$-partition}, we mean a collection 
$\cF= \{\alpha_1, \alpha_2,\cdots, \alpha_r \}$ of functions given by
$$
\alpha_i ~:~ A_i \longrightarrow G
$$
for $i\in \{1, 2, \cdots, r\}$, where $\{A_1, A_2, \cdots, A_r\}$ 
is a partial partition of $X$. Let $\tilde{Q}_n(G)$ be the 
set of all partial $G$-partitions of $X$. 
We have a map
\begin{align*}
 \psi : \tilde{Q}_n(G) &~\longrightarrow~ Q_n \\
       \{\alpha_1, \alpha_2,\cdots, \alpha_r \} 
       &~\longmapsto~ \{A_1, A_2, \cdots, A_r\}.
\end{align*}
 Let  $\cF = \{\alpha_1, \alpha_2,\cdots, \alpha_r\}$ 
 and $\cG= \{\beta_1, \beta_2, \cdots, \beta_s\}$ be 
 two partial $G$-partitions of $X$, 
 where $\alpha_i : A_i \rightarrow G$ and 
 $\beta_j : B_j \rightarrow G$ for $1 \leq i \leq r$ and $1 \leq j \leq s$. 
 Suppose that  $\psi(\cF) \leq \psi(\cG)$. 
 Then for each $j \in \{1,2, \cdots, s\}$, 
 there exists a subset $M_j \subseteq  \{1,2, \cdots, r\}$ such that
$$
B_j ~=~ \bigcup_{i \in M_j} A_i
$$ 
We say that $\beta_j$ is a  $G$-linear combination of $\alpha_i$'s if
$$
\beta_j|_{A_i} ~=~ \lambda_i \alpha_i
$$
for some scalars $\lambda_i \in G$, $i \in M_j$
and we write $\beta_j$ as
$$
\beta_j ~=~ \sum_{i \in M_j} \lambda_i \alpha_i.
$$

%Then we can write
%$$
%\beta_j = \sum_{i} \beta_j(x_i) e_i.
%$$

One defines a {\bf preorder} on $\tilde{Q}_n(G)$ as follows: we write 
$$
\cF ~\leq~ \cG,
$$ 
if $\psi(\cF) \leq \psi(\cG)$ and  $\beta_j$ is a $G$-linear combinations 
of $\alpha_i$'s for every $j \in \{1, 2, \cdots, s\}$. 
Clearly, the relation $\leq$ is reflexive and transitive, 
but not antisymmetric unless $G=\{1\}$. 
Hence we define an equivalence relation on $\tilde{Q}_n(G)$ by
$$
\cF ~\sim \cG \phantom{m}\text{if } \cF 
~\leq~ \cG \text{  and } \cG ~\leq~ \cF.
$$
The preorder $\leq$ on $\tilde{Q}_n(G)$ induces a partial order on the quotient set 
$$
Q_n(G) ~=~ \tilde{Q}_n(G)/\sim.
$$
In \cite{Do}, Dowling proved that $Q_n(G)$ is a geometric lattice. 
Note that if $G$ is the trivial group, then $Q_n(G) \cong Q_n \cong P_{n+1}$. 
The equivalence class of $\cF$ is usually denote by $(\cF)$. 

Let $E_i = \{x_i\}$ for $i \in \{1, 2, \cdots, n\}$ and $e_i : E_i \rightarrow G $ 
be the function defined by $e_i(x_i)=1$. Also let $\varepsilon = \{e_1, e_2, \cdots , e_n\}$, 
a partial $G$-partition of $X$. 
The zero element $\hat{0}$ of $Q_n(G)$ is $(\varepsilon)$ 
and the unit element $\hat{1}$ is the equivalence class of empty partial $G$-partition.
The rank function on $Q_n(G)$ is given by
$$
r(\cF) ~=~ n- |\cF|.
$$
Dowling  proved the following result regarding the structure of intervals 
in $Q_n(G)$ (see \cite[Theorem 2]{Do}).
\begin{thm}\label{thmstrIntQn}
Let $(\cF) \in Q_n(G)$ be of co-rank $r$ and let $(\cG) \in Q_n(G)$, where
$\cG ~=~ \{\beta_j : 1\leq j \leq s\}$ 
with $\beta_j : B_j \rightarrow G$ for $1 \leq j \leq s$. 
\begin{itemize}
\item[i)] We have
	$$
	[(\cF), \hat{1}] ~\cong~ Q_r(G).
	$$
	
\item[ii)] Let $B_0 = X \backslash \(\bigcup_{j=1}^{s} B_j\)$ 
and $n_j = |B_j|$ for $0 \leq j \leq s$. Then we have
	$$
	[\hat{0}, (\cG)] ~\cong~ 
	Q_{n_0}(G) \times P_{n_1} \times P_{n_2} \times \cdots \times P_{n_s}.
	$$
\end{itemize}	
\end{thm} 
We have the following theorem regarding convolution of Whitney numbers (see \cite[Theorem 6]{Do}).
\begin{thm}\label{thmshiftconvWk}
For non-negative integers $n, s$, we have 
$$
\sum_{r\geq 0} W_m(n,r)w_m(r,s) 
~=~ \delta(n, s) 
~=~ \sum_{r\geq 0} w_m(n,r)W_m(r,s).
$$
\end{thm}

\medspace

\section{Proofs of main theorems}

\subsection{Proof of \thmref{ThmU}}
Let $\{a_i\}_{i=0}^{n}$ be a unimodal sequence of 
non-negative real numbers such that
\begin{equation}\label{eq0}
\sum_{i=0}^{n} (-1)^i a_i ~=~ 0.
\end{equation}
We will prove \thmref{ThmU} by contradiction. 
Suppose there exists an integer $k \geq 0$ such that
\begin{equation}\label{eq1}
	\sum_{i=0}^{2k} (-1)^i a_i ~<~ 0.
\end{equation}
Since the sequence $\{a_i\}_{i=0}^{n}$ is unimodal, 
there exists an index $0 \leq j \leq n$ such that
$$
a_0 ~\leq~ a_1 ~\leq~ \cdots ~\leq~ a_{j-1} 
~\leq~ a_j ~\geq~ a_{j+1} ~\geq~ \cdots ~\geq~ a_n.
$$
If possible $2k \leq j$, then we have
$$
\sum_{i=0}^{2k} (-1)^i a_i 
~=~ a_0 \underbrace{-a_1 + a_2}_{\geq 0} - \cdots 
\underbrace{-a_{2k-1} + a_{2k}}_{\geq 0} 
~\geq~ 0.
$$
Hence we get $2k > j$. 
We split the following sum into two parts:
\begin{equation}\label{eqsplit}
\sum_{i=0}^{n} (-1)^i a_i 
~=~ \sum_{i=0}^{2k} (-1)^i a_i ~+~ 
\sum_{i=2k+1}^{n} (-1)^i a_i.
\end{equation} 
From the fact that $2k > j$, we get
\begin{equation}\label{eqsp2}
\sum_{i=2k+1}^{n} (-1)^i a_i ~=~
\underbrace{-a_{2k+1} + a_{2k+2}}_{\leq 0} 
- \cdots + (-1)^n a_n 
~\leq~ 0. 
\end{equation}
From \eqref{eq1}, \eqref{eqsplit} and \eqref{eqsp2}, we get 
$$
\sum_{i=0}^{n} (-1)^i a_i  ~<~  0,
$$
a contradiction to \eqref{eq0}. 
This completes the proof of first part of \thmref{ThmU}. 
The second part follows in a similar manner.
Suppose that for some positive integer $k$,
\begin{equation*}
	\sum_{i=0}^{2k-1} (-1)^i a_i  ~>~  0.
\end{equation*}
Then as before, we must have $2k-1 > j$. 
Now, we split the sum as
$$
\sum_{i=0}^{n} (-1)^i a_i 
~=~ \sum_{i=0}^{2k-1} (-1)^i a_i ~+~ 
\sum_{i=2k}^{n} (-1)^i a_i
$$
and note that
$$
\sum_{i=2k}^{n} (-1)^i a_i ~=~ \underbrace{a_{2k} -a_{2k+1}}_{\geq 0} 
+ \cdots +(-1)^{n} a_n 
~\geq~ 0.
$$
This leads to 
$$
\sum_{i=0}^{n} (-1)^i a_i ~>~ 0,
$$
a contradiction to \eqref{eq0}. 
This completes the proof of \thmref{ThmU}. \qed

\medspace

Now as an application of Heron-Rota-Welsh conjecture 
(now a theorem due to Adiprasito-Huh-Katz) and \thmref{ThmU}, 
we give a proof  Brun's inequality for a geometric lattice.

\subsection{Proof of \thmref{ThmB}}
The outline of proof is as follows : First we notice that 
it is equivalent to prove Brun's inequality for lattices 
associated to simple matroids. Then, we see that 
Whitney numbers of the first kind can be expressed 
in terms of values of the M\"obius function on such lattices. 
Finally, we complete the proof of \thmref{ThmB}
using the result of Adiprasito-Huh-Katz and \thmref{ThmU}.

\vspace{1mm}

Let $M = (E, \cF)$ be a matroid. The set $L(M)$ of flats of $M$ 
is a poset ordered by inclusion. 
One can show that $L(M)$ is a geometric lattice. 
Infact, every geometric lattice arises in this way 
upto isomorphism (see \cite[Theorem 1.7.5]{Ox}). 
More precisely, given any geometric lattice $L$, 
there exists a matroid $M$ such that 
$$
L ~\cong~ L(M).
$$
Further, one can suppose that $M$ is simple, 
since $L(M) \cong L(\hat{M})$.

Now onwards, let $M$ denote a simple matroid and 
$L(M)$ denote the lattice of flats of $M$. 
Let $\mu = \mu_{L(M)}$ denote the M\"obius function on $L(M)$. 
Unlike in the classical case where 
the M\"obius function takes the values $0$ or $\pm 1$, 
computing values of M\"obius functions on posets is a difficult problem. 
In this context, we have the following important formula 
(see \cite{Ha, Wei, Wh}, \cite[Prop. 7.1.4]{Za}):
\begin{equation}\label{valMob}
	\mu(\emptyset, F) ~=~ 
	\sum_{\substack{A \subseteq F\\ \overline{A} = F}} (-1)^{|A|}. 
\end{equation}
If we consider the polynomial
$$
\chi_{L(M)} (\lambda) ~=~ 
\sum_{ F \in L(M)} \mu(\emptyset, F) \lambda^{r(M) - r(F)}
$$
called the characteristic polynomial of the lattice $L(M)$, then one has
$$
\chi_{L(M)} (\lambda) ~=~ \chi_{M}(\lambda).
$$
This can be shown by using the formula \eqref{valMob} as follows:
\begin{equation}\label{eq2}
\begin{split}
\sum_{ F \in L(M)} \mu(\emptyset, F) \lambda^{r(M) - r(F)} 
~=~ \sum_{ F \in L(M)}\sum_{\substack{A \subseteq F\\ \overline{A} = F}} 
(-1)^{|A|} \lambda^{r(M) - r(F)}.
\end{split}
\end{equation}
Note that $r(F) = r(A)$, since $\overline{A} = F$. 
Hence the right hand side of \eqref{eq2} is equal to
$$
\sum_{ F \in L(M)}\sum_{\substack{A \subseteq F\\ \overline{A} = F}} 
(-1)^{|A|} \lambda^{r(M) - r(A)} 
~=~ 
\sum_{A \subseteq E} (-1)^{|A|} \lambda^{r(M) - r(A)}
~=~ \chi_{M}(\lambda).
$$
Thus we have
\begin{equation}\label{eqmuwi}
\sum_{ F \in L(M)} \mu(\emptyset, F) \lambda^{r(M) - r(F)} 
~=~ \sum_{i=0}^{r(M)} w_i(M) \lambda^{r(M) - i}.
\end{equation}
Now by comparing the coefficients on both sides, we get
\begin{equation}\label{eqmuw}
w_i(M) ~=~ 
\sum_{\substack{ F \in L(M)\\ r(F) = i}} \mu(\emptyset, F).	
\end{equation}
This in turn gives 
\begin{equation}\label{eqmuFk}
\sum_{\substack{ F \in L(M)\\ r(F) \leq k}} \mu(\emptyset, F) 
~=~
 \sum_{i=0}^{k} w_i(M).
\end{equation}
To complete the proof of \thmref{ThmB}, 
we need to show that $\sum_{i} w_i(M)$ is `alternating in sign' i.e.,
$$
(-1)^k  \sum_{i=0}^{k} w_i(M) ~\geq~ 0.
$$
Thus, it is important to know the 
`sign' of Whitney numbers of the first kind. 
In \cite[Theorem 4]{Ro64}, Rota proved the following fundamental result 
on the sign of M\"obius function of a geometric lattice : 
If $L$ is a geometric lattice, then the M\"obius function 
is non-zero and alternates in sign, 
i.e., for any $x \leq y$ in $L$,
\begin{equation}
	(-1)^{r(y)-r(x)} \mu_L(x, y) ~>~ 0.
\end{equation}
From this result and \eqref{eqmuw}, it follows that the numbers 
$w_i(M)$ alternates in sign i.e., 
\begin{equation}\label{eqsignW}
(-1)^i w_i(M) ~\geq~ 0.
\end{equation}
But this alone will not be sufficient for the completion of the proof, 
we need to establish that the sums
$$
\sum_{i} w_i(M)
$$
alternates in sign. By substituting $\lambda = 1$ in \eqref{eqmuwi}, we get
\begin{equation}\label{eqcond1}
\sum_{i=0}^{r(M)} w_i(M) 
~=~ \sum_{ F \in L(M)} \mu(\emptyset, F) 
~=~ 0.
\end{equation}
By using the fact that $|w_i(M)| = (-1)^i w_i(M)$, we rewrite the sum in \eqref{eqmuFk} as
$$
\sum_{\substack{ F \in L(M)\\ r(F) \leq k}} \mu(\emptyset, F) 
~=~
\sum_{i=0}^{k} w_i(M) 
~=~ \sum_{i=0}^{k} (-1)^i|w_i(M)|
$$
Now using the result of Adiprasito-Huh-Katz, \eqref{eqsignW} and  \eqref{eqcond1}, 
we apply \thmref{ThmU} to conclude that
 $$
 \sum_{\substack{F \in L(M) \\ r(F) \leq 2k}} \mu(\emptyset, F) 
 ~\geq~ 0
 \phantom{mm} \text{and} \phantom{mm}
 \sum_{\substack{F \in L(M) \\ r(F) \leq 2k-1}} \mu(\emptyset, F) 
 ~\leq~ 0.
 $$
This completes the proof of \thmref{ThmB}. \qed

\medspace

\section{A sieve on a geometric lattice}
In this section, we will establish a sieve for a geometric lattice.

\subsection{Proof of \thmref{thmSift}}
Let $L$ be a geometric lattice of rank $n$ 
and $\cA \subseteq L$. Also let $\cT$ be a set of atoms in $L$ 
and $\tau$ be the join of atoms in $\cT$.
Set
$$
\cS\(\cA, \cT\) ~=~ \# \{a \in \cA ~:~ a \wedge \tau ~=~ \hat{0}\}.
$$
Then we have
\begin{equation}\label{eqSiftAy}
\begin{split}
\cS\(\cA, \cT\) 
&~=~ 
\sum_{a \in \cA \atop a \wedge \tau = \hat{0}} 1
~=~ \sum_{a \in \cA} \sum_{y \leq a \wedge \tau} \mu(\hat{0}, y) \\
&~=~ \sum_{y \leq \tau}  \mu(\hat{0}, y) \sum_{ a \in \cA \atop a \geq y} 1
~=~  \sum_{y \leq \tau}  \mu(\hat{0}, y) \#\cA_y,
\end{split}
\end{equation}
where $\cA_y = \{a \in \cA ~:~ a \geq y\}$. 
Suppose there exists a function $f: \N \cup \{0\} \rightarrow [0, \infty)$ such that
\begin{equation*}
\#\cA_y 
~=~ f(\text{cr}(y)) X ~+~ \cE(y)
\phantom{mm}\text{and}\phantom{mm}
\cE(y) ~\ll~ \text{cr}(y) f(\text{cr}(y)).
\end{equation*}
Then, we get
\begin{equation*}
 \sum_{y \leq \tau}  \mu(\hat{0}, y) \#\cA_y
 ~=~ X \sum_{y \leq \tau}  \mu(\hat{0}, y)  f(\text{cr}(y))  
 ~+~ \sum_{y \leq \tau}  \mu(\hat{0}, y) \cE(y).
\end{equation*}
Note that $y \leq \tau$ implies that $r(y) \leq r(\tau)$. 
By collecting the terms $y \leq \tau$ with same rank, we get
\begin{equation*}
\begin{split}
\sum_{y \leq \tau}  \mu(\hat{0}, y)  f(\text{cr}(y))
&~=~  \sum_{k=0}^{r(\tau)} \sum_{y \leq \tau \atop r(y)=k} \mu(\hat{0}, y) f(n-k) \\
&~=~  \sum_{k=0}^{r(\tau)} f(n-k) \sum_{y \in [\hat{0}, \tau]  \atop r(y)=k} \mu(\hat{0}, y) \\
&~=~  \sum_{k=0}^{r(\tau)} f(n-k) w_k\([\hat{0}, \tau]\),
\end{split}
\end{equation*}
where $w_k\([\hat{0}, \tau]\)$ denotes the $k$-th Whitney number of the first kind 
of the geometric lattice $[\hat{0}, \tau]$ and 
note that any interval in a geometric lattice is geometric.
By arguing in a similar way and using the fact that $(-1)^{r(y)} \mu(\hat{0}, y) \geq 0$, 
we deduce that
\begin{equation*}
\sum_{y \leq \tau}  \mu(\hat{0}, y) \cE(y) ~\ll~
\sum_{k=0}^{r(\tau)} (n-k)f(n-k) |w_{k}([\hat{0}, \tau])|.
\end{equation*}
This completes the proof of \thmref{thmSift}. \qed

\medspace

\section{Sifted sets and Shifted convolution of Whitney numbers}
In this section, we study relations between the cardinality of sifted sets 
and shifted convolution of Whitney numbers in the special case of Dowling lattices.
Consider $L=\cA= Q_n(G)$ and $\cT$ be a set of atoms whose join is $\tau$. 
Also let $k = r(\tau)$. Then from \thmref{thmstrIntQn}, we have
\begin{equation}
\#\cA_y ~=~ \#[y, \hat{1}] ~=~ \#Q_{n-r(y)}(G)
\end{equation}
We set 
$$
X ~=~ \#Q_{n}(G)
\phantom{m}\text{and}\phantom{m}
f(m) ~=~ \frac{\#Q_{m}(G)}{\#Q_{n}(G)} 
\phantom{m}\text{for}\phantom{m}
m  \in \N \cup \{0\}.
$$
From \eqref{eqSiftAy}, we have
\begin{equation*}
\begin{split}
\cS\(\cA, \cT\) 
&~=~ \#Q_{n}(G) \sum_{s=0}^{k} f(n-s) 
\sum_{y \leq \tau \atop r(y)=s } \mu(\hat{0}, y) 
~=~ \sum_{s=0}^{k} \#Q_{n-s}(G) 
\sum_{y \leq \tau \atop r(y)=s } \mu(\hat{0}, y) \\
&~=~ \sum_{s=0}^{k} \#Q_{n-k+s}(G) 
\sum_{y \leq \tau \atop r(y)=k-s } \mu(\hat{0}, y).
\end{split}
\end{equation*}
For simplicity, we suppose that all blocks of $\tau$ are trivial except the zero-block, 
so that by \thmref{thmstrIntQn}, we have
$$
[\hat{0}, \tau]  ~\cong~ Q_k(G)
\phantom{mm}\text{and}\phantom{mm} 
\sum_{y \leq \tau \atop r(y)=k-s } \mu(\hat{0}, y) ~=~ w_m(k,s).
$$
Thus we get
\begin{equation}\label{eqSiftShift}
\begin{split}
\cS\(\cA, \cT\)
~=~ \sum_{s=0}^{k} \#Q_{n-k+s}(G) \cdot w_m(k,s) 
~=~ \sum_{s=0}^{k} \sum_{r=0}^{n-k+s}  w_m(k, s) W_m(n-k+s, r).
\end{split}
\end{equation}
Hence the estimation of cardinality of sifted sets leads us 
to the study of {\bf shifted convolution}
of Whitney numbers of Dowling lattices:
$$
\sum_{k \geq 0} w_m(n, k) W_m(k+s,t).
$$

\medspace

\subsection{Proof of \propref{propshiftW}}
Let $n , s,t$ be non-negative integers and set
\begin{equation}\label{eqshiftconv}
c_{n,t}(s) ~=~ \sum_{k \geq 0} w_m(n, k) W_m(k+s,t).
\end{equation}
Notice that the sum in \eqref{eqshiftconv} is a finite sum, 
since $w_m(n, k) = 0$ whenever $k > n$.
We know that (see \cite[Theorem 7]{Do}) 
the numbers $w(n,k)$ satisfy the recurrence relation
$$
w_m(n,k) ~=~ w_m(n-1, k-1) ~-~ (1+m(n-1)) w_m(n-1, k)
$$
and hence we get
\begin{equation}\label{eqreccnr}
c_{n,t}(s) ~=~ c_{n+1, t}(s-1) ~+~ (1+mn) c_{n,t}(s-1), 
\phantom{mm} s \in \N.
\end{equation}
Consider the rational generating function
$$
F_{n,t}(x) ~=~ \sum_{s=0}^{\infty} c_{n,t}(s) x^s.
$$
From \eqref{eqreccnr} and \thmref{thmshiftconvWk},
we get
\begin{equation}
F_{n,t}(x) 
~=~ \frac{1-(1+m(n-1))x}{x} \cdot F_{n-1,t}(x) 
~-~ \frac{\delta(n-1, t)}{x}.
\end{equation}
By induction on $n$, we deduce that
\begin{equation}\label{eqFnrind}
\begin{split}
F_{n,t}(x) ~=~ \prod_{i=0}^{n-1} \frac{1-(1+im)x}{x} \cdot F_{0,t}(x) 
~-~ \sum_{j=0}^{n-1} \prod_{i=j+1}^{n-1} \frac{1-(1+im)x}{x} \cdot \frac{\delta(j,t)}{x}.
\end{split}
\end{equation}
From \cite[Theorem 5]{Be}, we get
\begin{equation}\label{eqF0recWm}
F_{0,t}(x) ~=~ \sum_{s=0}^{\infty} c_{0,t}(s) x^s 
~=~ \sum_{s=t}^{\infty} W_m(s,t) x^s 
~=~ \frac{1}{x} \prod_{i=0}^{t} \frac{x}{1-(1+im)x}.
\end{equation}
From \eqref{eqFnrind} and \eqref{eqF0recWm}, we have
\begin{equation}\label{eqFnr}
\begin{split}
F_{n,t}(x) ~=~ 
\frac{1}{x} \prod_{i=0}^{n-1} \frac{1-(1+im)x}{x} \cdot
 \prod_{j=0}^{t} \frac{x}{1-(1+ jm)x}
~-~ \sum_{j=0}^{n-1} \prod_{i=j+1}^{n-1} \frac{1-(1+im)x}{x} \cdot \frac{\delta(j,t)}{x}.
\end{split}
\end{equation}
If $t < n$, then from \eqref{eqFnr}, we get
\begin{equation}
F_{n,t}(x) ~=~ 
\frac{1}{x} \prod_{i=t+1}^{n-1} \frac{1-(1+im)x}{x}
~-~ \prod_{i=t+1}^{n-1} \frac{1-(1+im)x}{x} \cdot \frac{1}{x}
~=~ 0.
\end{equation}
Thus we conclude that if $t < n$, then we have
\begin{equation}\label{eqcnr<n}
c_{n,t}(s) ~=~ 0 
\phantom{m}\forall~~ s \in \N \cup\{0\}.
\end{equation}
%If $r=n$, then from \eqref{eqFnr}, we get
%\begin{equation}
%F_{n,r}(x) ~=~ \frac{1}{1-(1+mn)x} ~=~ \sum_{N=0}^{\infty} (1+mn)^N x^N.
%\end{equation}
%Hence if $r=n$, then we get
%\begin{equation}
%c_{n,r}(N) ~=~ (1+mn)^N
%\phantom{m}\text{for}\phantom{m} N \geq 0.
%\end{equation}
Suppose that $t\geq n$. Then from \eqref{eqFnr}, we get
\begin{equation}\label{eqFnr>n}
F_{n,t}(x) ~=~ \frac{1}{x} \prod_{j=n}^{t} \frac{x}{1-(1+jm)x}.
\end{equation}
This completes the proof of \propref{propshiftW}. \qed

\medspace

\section{Sifted sets and Dowling numbers}\label{SDmr}
In this section, we study relations between certain sifted sets 
and generalized Dowling numbers. Further, 
we derive an asymptotic expression for the generalized Dowling numbers.

For any non-negative integer $r$, one defines (see \cite{CJ, Me}) 
$r$-Whitney numbers of the second kind of $Q_n(G)$ by
$$
(mx+r)^n ~=~ \sum_{k=0}^{n} m^k W_{m,r}(n, k) (x)_k,
$$
where $(x)_k = x(x-1) \cdots (x-k+1)$ denotes the $k$-th falling factorial.
Then the numbers $W_{m,r}(n, k)$ satisfy the recurrence relation
$$
W_{m,r}(n, k) ~=~ W_{m,r}(n-1, k-1) + (km+r) W_{m,r}(n-1, k)
$$
and one can show that its rational generating function is given by 
\begin{equation}\label{eqWmrgen}
\sum_{s=0}^{\infty} W_{m,r}(s, k) x^s 
~=~ \frac{1}{x} \prod_{i=0}^{k} \frac{x}{1-(r+im)x}.
\end{equation}
From \eqref{eqcnr<n}, \eqref{eqFnr>n} and \eqref{eqWmrgen}, 
we conclude that
\begin{equation}\label{eqcnrWm}
c_{n,t}(s) ~=~ W_{m, 1+mn}(s, t-n)
\phantom{m}\forall~~ s \in \N \cup\{0\}.
\end{equation}
From \eqref{eqSiftShift} and \eqref{eqcnrWm}, we get
\begin{equation}\label{eqsiftr-Dn}
\begin{split}
\cS(\cA, \cT) 
&~=~ \sum_{r} \sum_{s} w_{m}(k,s) W_m(n-k+s, r) 
~=~ \sum_{r} c_{k,r}(n-k) \\
&~=~ \sum_{r} W_{m, 1+mk}(n-k, r-k)
~=~ \sum_{r=0}^{n-k}  W_{m, 1+mk}(n-k, r) .
\end{split}
\end{equation}
The Dowling numbers associated with $Q_n(G)$ 
(see \cite[p. 22]{Be} \cite[Sec. 5]{CJ}) are defined by 
$$
D_m(n) ~=~ \sum_{k=0}^{n} W_m(n,k).
$$
More generally, one defines $r$-Dowling numbers associated with $Q_n(G)$ by
$$
D_{m, r}(n) ~=~ \sum_{k=0}^{n} W_{m, r}(n,k).
$$
From \eqref{eqsiftr-Dn}, we conclude that
\begin{equation}
\cS(\cA, \cT) ~=~ D_{m, 1+mk}(n-k).
\end{equation}

In the next subsection, we will give a proof of \thmref{thmexpDn} 
which gives an explicit expression for $r$-Dowling numbers $D_{m, r}(n)$ 
using the Cauchy's residue theorem and the method of steepest descent/saddle-point method 
(see for example \cite[Sec. 4.2]{Be} and \cite{BS}) and 
the expression obtained is uniform in $m$ and $r$.

\subsection{Proof of \thmref{thmexpDn}}
The exponential generating function of $D_{m,r}(n)$ 
(see \cite[Eq. 26]{CJ}) is given by
\begin{equation*}
\sum_{n=0}^{\infty} D_{m,r}(n) \frac{z^n}{n!} 
~=~ \exp\( rz + \frac{e^{mz} -1}{m}\)
\phantom{m}\text{for}\phantom{m} z \in \C.
\end{equation*}
By Cauchy's residue theorem, we get
\begin{equation}\label{eqDmrint}
\begin{split}
D_{m,r}(n) 
&~=~ \frac{n!}{2\pi i} \int_{|z|= \delta} 
\frac{\exp\(rz+ \frac{e^{mz} -1}{m}\)}{z^{n+1}} ~dz \\
&~=~ \frac{n!}{2 \pi \delta^n} \int_{-\pi}^{\pi} 
\exp\( r \delta e^{i \theta} - i n \theta + \frac{e^{m \delta e^{i \theta}} - 1}{m}\) ~ d\theta,
\end{split}
\end{equation}
where $\delta > 0$ is a constant which will be chosen later. Set
\begin{equation*}
G(z) ~=~ r \delta e^{i z} - i n z + \frac{e^{m \delta e^{i z}} - 1}{m}.
\end{equation*}
Then $G$ is an entire function and the power series expansion of $G$ 
around $z=0$ is given by
$$
G(z) ~=~ G(0) + G'(0) z + G''(0) \frac{z^2}{2} + H(z) z^3
$$
for some entire function $H$. We choose $\delta$ such that $G'(0)=0$ i.e., 
$\delta$ is a positive constant such that
\begin{equation}\label{eqdeltaG'0}
r \delta ~+~ \delta e^{m \delta}  ~=~ n.
\end{equation}
Then we have
\begin{equation}\label{eqg012}
g_0 ~:=~ G(0) ~=~ r \delta + \frac{e^{m\delta} -1}{m}
%g_1 ~:=~ G'(0) ~=~ r \delta -n + \delta e^{m \delta} 
\phantom{m}\text{and}\phantom{m}
g_2 := -\frac{G''(0)}{2} ~=~  \frac{n+ m \delta^2 e^{m \delta}}{2}.
\end{equation}
Let
\begin{equation}\label{eqchoeps}
	\e ~=~ \frac{(\log n)^5}{n^{1/2}}.
\end{equation}
We split the integral in \eqref{eqDmrint} into three parts as follows:
\begin{equation}\label{eqDmrint3parts}
D_{m, r}(n) ~=~ \frac{n!}{2 \pi \delta^n} 
\( \int_{-\pi}^{-\epsilon} + \int_{-\epsilon}^{\epsilon} + \int_{\epsilon}^{\pi} \)
\exp\(g_0 - g_2 \theta^2 + H(\theta) \theta^3\) ~ d\theta.
\end{equation}
Now we will show that
\begin{equation}\label{eqintpipie}
\( \int_{-\pi}^{-\epsilon} + \int_{\epsilon}^{\pi} \) 
\exp\( - g_2 \theta^2 + \Re\{H(\theta)\} \theta^3\) ~ d\theta 
~\ll~ \exp\(-\frac{(\log n)^9}{20}\).
\end{equation}
Note that 
\begin{equation}\label{eq-g2+H}
\begin{split}
- g_2 \theta^2 + \Re\{H(\theta)\} \theta^3 
&~=~ \Re\{G(\theta)\} - \Re\{G(0)\} \\
&~=~ r \delta ( \cos \theta - 1) 
~+~ \frac{e^{m \delta \cos \theta} \cos\( m \delta \sin \theta\) - e^{m \delta}}{m}.
\end{split}
\end{equation}
Note that
\begin{equation}\label{eqcos-1}
\cos \theta -1 ~\leq~ - \frac{\theta^2}{5}
\phantom{m}\text{for}\phantom{m} \theta \in [-\pi, \pi]
\end{equation}
and
\begin{equation}\label{eqemdelcos}
\begin{split}
e^{m \delta \cos \theta} \cos\( m \delta \sin \theta\) - e^{m \delta}
&~\leq~ e^{m \delta \cos \theta} - e^{m \delta}  
~=~ \sum_{k=1}^{\infty} \frac{(m \delta \cos \theta)^k - (m\delta)^k}{k!} \\
&~=~ \sum_{k=1}^{\infty} \frac{(m \delta)^k}{k!} \(\cos^k \theta  - 1\).
\end{split}
\end{equation} 
Thus for  $\theta \in [0, \pi/2]$, we have
\begin{equation}\label{eqemdelx<pi/2}
 e^{m \delta \cos \theta} - e^{m \delta} 
 ~\leq~ 
 (\cos \theta -1) \sum_{k=1}^{\infty} \frac{(m \delta)^k}{k!} 
 ~\leq~ -\frac{\theta^2}{5} \(e^{m \delta}-1\).
\end{equation}
If $\theta \in [\pi/2, \pi]$, then we have
\begin{equation}\label{eqemdelx>pi/2}
e^{m \delta \cos \theta} - e^{m \delta} 
~\leq~ 1 - e^{m \delta}.
\end{equation}
From \eqref{eqemdelcos}, \eqref{eqemdelx<pi/2} and \eqref{eqemdelx>pi/2}, 
we deduce that
\begin{equation}\label{eqemdelcos-emd}
e^{m \delta \cos \theta} \cos\( m \delta \sin \theta\) - e^{m \delta}
~\leq~ 1- \frac{\theta^2}{10} \(e^{m \delta} -1\)
\phantom{m}\text{for}\phantom{m} \theta \in [-\pi, \pi].
\end{equation}
From \eqref{eq-g2+H}, \eqref{eqcos-1} and \eqref{eqemdelcos-emd}, we get
\begin{equation}
- g_2 \theta^2 + \Re\{H(\theta)\} \theta^3 
~\leq~  
-  \frac{r \delta \theta^2}{5} + \frac{1}{m} - \frac{\theta^2}{10m} \(e^{m \delta}-1\).
\end{equation}
Thus we have
\begin{equation}\label{eqint-pipirdel}
\( \int_{-\pi}^{-\epsilon} + \int_{\epsilon}^{\pi} \) 
\exp\( - g_2 \theta^2 + \Re\{H(\theta)\} \theta^3\) ~ d\theta
~\leq~ 2 \int_{\e}^{\pi} \exp\( 1 + \frac{\theta^2}{10} - \(2 r \delta + \frac{e^{m \delta}}{m} \) \frac{\theta^2}{10}  \) ~ d\theta .
\end{equation}
Suppose that $n > e^m$. From \eqref{eqdeltaG'0}, we get
$$
\delta ~\leq~ \frac{\log n}{m}.
$$
%From \eqref{eqdeltaG'0}, we also note that if $r \delta < n/2$, then
%$$
%\frac{n}{2} 
%~\leq~ \delta e^{m \delta} 
%~\leq~ \frac{\log n}{m} e^{m \delta} 
%$$
%and this implies that
%$$
%e^{m \delta} ~\geq~ \frac{m}{2} \frac{n}{\log n}
%$$
%provided $r \delta < n/2$.
Further, we get 
\begin{equation}
2 r \delta + \frac{e^{m \delta}}{m} ~\geq~ \frac{n}{2 \log n}.
\end{equation}
Therefore, we have
\begin{equation}\label{eqintpie}
\begin{split}
\int_{\e}^{\pi} \exp\( - \(2 r \delta + \frac{e^{m \delta}}{m} \) \frac{\theta^2}{10}  \) ~ d\theta
&~\leq~ 
\int_{\e}^{\pi} \exp\( - \frac{n}{ 20\log n} \theta^2 \) ~ d\theta \\
&~\leq~
\pi \exp\( - \frac{n}{ 20\log n} \e^2 \)
~\leq~ 
\pi \exp\( - \frac{(\log n)^9}{20}\).
\end{split}
\end{equation}
Now the estimate in \eqref{eqintpipie} follows from \eqref{eqint-pipirdel} and \eqref{eqintpie}. 
Thus, we get
\begin{equation}\label{eqDmreps}
D_{m,r}(n) ~=~ \frac{e^{g_0} n!}{2 \pi \delta^n} 
\(\int_{-\e}^{\e} \exp \(-g_2 \theta^2 + H(\theta) \theta^3 \) ~d\theta 
~+~ O\(\exp\( - \frac{(\log n)^9}{20}\) \)\).
\end{equation}
Consider the entire function
$$
h(z) ~=~ \exp\(H(z) z^3\)
\phantom{m}\text{for}\phantom{m} z \in \C.
$$
and write
$$
h(z) ~=~ 1 + \tilde{h}(z) z
$$
for some entire function $\tilde{h}$. Thus we have
\begin{equation}\label{eqintHtildeh}
\begin{split}
\int_{-\e}^{\e} \exp \(-g_2 \theta^2 + H(\theta) \theta^3 \) ~d\theta
&~=~ \int_{-\e}^{\e} \exp \(-g_2 \theta^2 \) 
\(1+ \tilde{h}(\theta) \theta\) ~d\theta \\
&~=~ \int_{-\e}^{\e} \exp \(-g_2 \theta^2 \)  ~d\theta 
~+~ \int_{-\e}^{\e} \exp \(-g_2 \theta^2 \) \tilde{h}(\theta) \theta ~d\theta .
\end{split}
\end{equation}
By noting that $|e^w -1| \leq \max\{|w|, |w|e^{\Re\{w\}} \} \leq |w| e^{|w|} $ 
for any $w \in \C$, we get
\begin{equation}\label{eqhtildeupp}
|\tilde{h}(z)| 
~=~ \left| \frac{h(z) -1}{z} \right|
~=~ \left|\frac{e^{H(z) z^3} -1}{z} \right|
~\leq~ |H(z)| |z|^2 e^{|H(z)||z|^3} 
\phantom{m}\text{for}\phantom{m} z \in \C\backslash\{0\}
\end{equation}
and by continuity, the inequality in \eqref{eqhtildeupp} holds for all $z \in \C$.
For $z \in \C\backslash\{0\}$,  we have
\begin{equation*}
\begin{split}
H(z) &~=~ \frac{G(z) - g_0 +g_2 z^2}{z^3} \\
&~=~ 
\frac{n}{z^3} \sum_{k=3}^{\infty} \frac{(iz)^k}{k!} ~+~
\frac{i m \delta^2 e^{m \delta}}{z^2} \sum_{k=2}^{\infty} \frac{(iz)^k}{k!}
~+~ \frac{m \delta^2 e^{m \delta}}{2 z^3}  
\( \sum_{k=2}^{\infty} \frac{(iz)^k}{k!} \)^2
~+~ \frac{e^{m \delta}}{ m z^3} 
\sum_{k=3}^{\infty} \frac{(m \delta (e^{iz}-1))^k}{k!}
\end{split}
\end{equation*}
Thus for $\theta \in [-\e, \e] \backslash\{0\}$, we get
\begin{equation}
\begin{split}
|H(\theta)| 
&~\ll~ n
~+~ m \delta^2 e^{m \delta} 
~+~ \frac{e^{m \delta}}{ m |\theta|^3} 
\sum_{k=3}^{\infty} \frac{(m \delta |\theta|)^k}{k!} \\
&~\ll~ n
~+~ m \delta n
~+~ m^2 \delta^3 e^{(1+\e) m \delta}.
\end{split}
\end{equation}
From \eqref{eqdeltaG'0} and \eqref{eqchoeps}, we get
\begin{equation}\label{eqHuppeps}
|H(\theta)| ~\ll~ (1+ m\delta) n ~+~ m^2 \delta^2 n 
 ~\ll~ n (\log n)^2,
\end{equation}
where the implied constant is absolute.
From \eqref{eqhtildeupp} and \eqref{eqHuppeps}, we get
\begin{equation}
|\tilde{h}(\theta)| ~\ll~ (\log n)^{12} 
\phantom{m}\text{for}\phantom{m} 
\theta \in [-\e, \e].
\end{equation}
Hence we get
\begin{equation}\label{eqintg2tildeh}
\int_{-\e}^{\e} \exp \(-g_2 \theta^2 \) |\tilde{h}(\theta) \theta| ~d\theta 
~\ll~ \frac{(\log n)^{12}}{g_2}.
\end{equation}
Also, we have
\begin{equation}\label{eqintg2e}
\begin{split}
\int_{-\e}^{\e} \exp \(-g_2 \theta^2 \)  ~d\theta 
&~=~ \frac{1}{\sqrt{g_2}} \int_{-\e \sqrt{g_2}}^{\e \sqrt{g_2}}  e^{-u^2} ~du 
~=~  \frac{1}{\sqrt{g_2}} \int_{-\infty}^{\infty} e^{-u^2} ~du 
~-~ \frac{2}{\sqrt{g_2}} \int_{\e \sqrt{g_2}}^{\infty} e^{-u^2} ~du \\
&~=~ \frac{\sqrt{\pi}}{\sqrt{g_2}} 
~+~ O\( \frac{e^{- \e^2 g_2 }}{\e g_2}\).
\end{split}
\end{equation}
From \eqref{eqDmreps}, \eqref{eqintHtildeh}, \eqref{eqintg2tildeh} and \eqref{eqintg2e}, we get
$$
D_{m,r}(n) ~=~ \frac{e^{g_0} n!}{\sqrt{4 \pi g_2} \delta^n}
\(1 ~+~ O\( \frac{e^{- \e^2 g_2 }}{\e \sqrt{g_2}} 
~+~ \frac{(\log n)^{12}}{\sqrt{g_2}}
~+~ \sqrt{g_2} \exp\( -\frac{(\log n)^9}{20}\) \)\).
$$
From \eqref{eqg012}, we note that $n/2 \leq g_2 \leq n \log n$. Hence for $n > e^m$, we have
\begin{equation*}
D_{m,r}(n) ~=~ \frac{e^{g_0} n!}{\sqrt{4 \pi g_2} \delta^n}
\(1 ~+~ O\(\frac{(\log n)^{12}}{\sqrt{n}}\)\),
\end{equation*}
where the implied constant is absolute. 
This completes the proof of \thmref{thmexpDn}. \qed 

\section{Concluding remarks}

We envisage a larger theory to emerge from these modest beginnings.
There are many questions in combinatorics that can be
formulated in some sort of sieve theoretic terms.  Coloring
problems come to mind where these ideas could have potential applications.  This is the hope for the future.

\medspace

\section*{Acknowledgments}
The authors would like to acknowledge Queen’s University, Canada for providing excellent atmosphere to work.

\medspace

\end{document}